\documentclass[a4paper,12pt]{article}
\usepackage[english]{babel}
\usepackage[latin1]{inputenc}
\usepackage[T1]{fontenc}
\usepackage{amsmath}
\usepackage{amsfonts}
\usepackage{amsthm}
\usepackage{amssymb}
\usepackage{dsfont}

\usepackage{epsfig}
\psfigdriver{dvips}
\usepackage{colt_court}
\begin{document}

\title{Model selection by resampling penalization}
\author{
Sylvain Arlot\\
Univ Paris-Sud, Laboratoire de Math\'ematiques d'Orsay, \\
Orsay Cedex, F-91405; CNRS, Orsay cedex, F-91405 \\
\texttt{sylvain.arlot@math.u-psud.fr}
\\
INRIA Futurs, Projet Select
}
\date{January 19, 2007}

\maketitle

\begin{abstract}
We present a new family of model selection algorithms based on the resampling heuristics. It can be used in several frameworks, do not require any knowledge about the unknown law of the data, and may be seen as a generalization of local Rademacher complexities and $V$-fold cross-validation. In the case example of least-square regression on histograms, we prove oracle inequalities, and that these algorithms are naturally adaptive to both the smoothness of the regression function and the variability of the noise level. Then, interpretating $V$-fold cross-validation in terms of penalization, we enlighten the question of choosing $V$. Finally, a simulation study illustrates the strength of resampling penalization algorithms against some classical ones, in particular with heteroscedastic data.
\end{abstract}

\section{Introduction}
Choosing between the outputs of many learning algorithms, from the prediction viewpoint, remains to estimate their generalization abilities. A classical method for this is penalization, that comes from model selection theory. Basically, it states that a good choice can be made by minimizing the sum of the empirical risk (how does the algorithm fits the data) and some complexity measure of the algorithm (called the penalty). 
The ideal penalty for prediction is of course the difference between the true and empirical risks of the output, but it is unknown in general. It is thus crucial to obtain tight estimates of such a quantity. 

Many penalties or complexity measures have been proposed, both in the classification and regression frameworks. Consider for instance regression and least-square estimators on finite-dimensional vector spaces (the models). When the design is fixed and the noise-level constant equal to $\sigma$, Mallows' $C_p$  penalty \cite{Mal:1973} (equal to $2 n^{-1} \sigma^2 D$ for a $D$-dimensional space, and it can be modified according to the number of models \cite{Bir_Mas:2002,Sau:2006}) has some optimality properties \cite{Shi:1981,KCLi:1987,Bar:2002}. However, such a penalty linear in the dimension may be terrible in an heteroscedastic framework (as shown by \eqref{eq.Em_penid} and experiment HSd2 in Sect.~\ref{sec.simus}).

In classification, the VC-dimension has the drawback of being independent of the underlying measure, so that it is adapted to the worst case. It has been improved with data-dependent complexity estimates, such as Rademacher complexities \cite{Kol:2001,Bar_Bou_Lug:2002} (generalized by Fromont with resampling ideas \cite{Fro:2004}), but they may be too large because they are still global complexity measures. The localization idea then led to local Rademacher complexities \cite{Bar_Bou_Men:2005,Kol:2006} which are tight estimates of the ideal penalty, but involve unknown constants and may be very difficult to compute in practice. On the other hand, the $V$-fold cross-validation (VFCV) is very popular for such purposes, but it is still poorly understood from the non-asymptotic viewpoint.

In this article, we propose a new family of penalties, based on Efron's bootstrap heuristics \cite{Efr:1979} (and its generalization to weighted bootstrap, \latin{i.e.} resampling). It is a localized version of Fromont's penalties, which does not involve any unknown constant, and is easy to compute (at the price of some loss in accuracy) in its $V$-fold cross-validation version. We define it in a much general framework, so that it has a wide range of application. As a first theoretical step, we prove the efficiency of these algorithms in the case example of least-square regression on histograms, under reasonable assumptions. Indeed, they satisfy oracle inequalities with constant almost one, asymptotic optimality and adaptivity to the regularity of the regression function. This comes from explicit computations that allow us to deeply understand why these penalties are working well. Then, we compare the ``classical'' VFCV with the $V$-fold penalties, enlightening how $V$ should be chosen. Finally, we illustrate these results with a few simulation experiments. In particular, we show that resampling penalties are competitive with classical methods for ``easy problems'', and may be much better for some harder ones (\latin{e.g.} with a variable noise-level).

\section{A general model selection algorithm} \label{sec.proc}
We consider the following general setting : $\X \times \Y$ is a measurable space, $P$ an unknown probability measure on it and $(X_1,Y_1), \ldots, (X_n,Y_n) \in \X \times \Y$ some data of common law $P$. Let $S$ be the set of predictors (measurable functions $\X \mapsto \Y$) and $\gamma : S \times (\X \times \Y) \mapsto \R$ a contrast function. Given a family $(\ERM_m)_{\mM_n}$ of data-dependent predictors, our goal is to find the one minimizing the prediction loss $P\gamma(t)$. We will extensively use this functional notation $Q \gamma(t) \egaldef \E_{(X,Y)\sim Q}[\gamma(t,(X,Y))]$, for any probability measure $Q$ on $\X \times \Y$. Notice that the expectation here is only taken w.r.t. $(X,Y)$, so that $Q \gamma(t)$ is random when $t=\ERM_m$ is random. Assuming that there exists a minimizer $\bayes \in S$ of the loss (the Bayes predictor), we will often consider the excess loss $\perte{t} = P\gamma(t) - P\gamma(\bayes) \geq 0$ instead of the loss.

Assume that each predictor $\ERM_m$ may be written as a function $\ERM_m(P_n)$ of the empirical distribution of the data $P_n = n^{-1} \sum_{i=1}^n \delta_{(X_i,Y_i)}$. 
The ideal choice for $\mh$ is the one which minimizes over $\M_n$ the true prediction risk $P \gamma(\ERM_m(P_n)) = P_n \gamma(\ERM_m(P_n)) + \penid(m)$ where the ideal penalty is equal to 
\begin{equation*} 
 \penid(m) = (P-P_n) \gamma(\ERM_m(P_n)) \enspace . \end{equation*}
The {\em resampling heuristics} (introduced by Efron \cite{Efr:1979}) states that the expectation of any functional $F(P,P_n)$ is close to its resampling counterpart $\Es F(P_n,\Pnb)$, where $\Pnb = n^{-1} \sum_{i=1}^n W_i \delta_{(X_i,Y_i)}$ is the empirical distribution $P_n$ weighted by an independent random vector $W \in [0;+\infty)^n$, with $\sum_i \E[W_i]=n$. The expectation $\Es[\cdot]$ means that we only integrate w.r.t. the weights $W$. 

We suggest here to use this heuristics for estimating $\penid(m)$, and plug it into the penalized criterion $P_n \gamma(\ERM_m) + \pen(m)$. This defines $\mh \in \M_n$ as follows. 
\begin{proc}[Resampling penalization] \label{def.proc.gal}
\begin{enumerate}
\item Choose a resampling scheme, \latin{i.e.} the law of a weight vector $W$.
\item Choose a constant $C \geq C_W \approx \left( n^{-1} \sum_{i=1}^n \E\paren{W_i-1}^2  \right)^{-1}$.
\item Compute the following resampling penalty for each $\mM_n$ :
\begin{equation*} 
\pen(m) = C \Es \croch{  P_n \gamma \paren{ \ERM_m \paren{ \Pnb} } - \Pnb \gamma \paren{\ERM_m \paren{\Pnb}} } \enspace .
\end{equation*}
\item Minimize the penalized empirical criterion to choose $\mh$ and thus $\ERM_{\mh}$ :
\begin{equation*} 
\mh \in \arg\min_{\mM_n} \left\{ P_n \gamma(\ERM_m(P_n)) + \pen(m) \right\} \enspace .
\end{equation*}
\end{enumerate}
\end{proc}

\begin{remark}
\begin{enumerate}
\item There is a constant $C \neq 1$ in front of the penalty, although there isn't any in Efron's heuristics, because we did not normalize $W$. The asymptotical value of the right normalizing constant $C_W$ may be derived from Theorem 3.6.13 in \cite{vdV_Wel:1996}. In the case example of histograms, we give a non-asymptotic expression for it \eqref{def.CW}. In general, we suggest to use some data-driven method to choose $C$ (see algorithm~\ref{def.proc.gal.pente}), whereas the resampling penalty only estimates the shape of the ideal one.
\item We allowed $C$ to be larger than $C_W$ because overpenalizing may be fruitful in a non-asymptotic viewpoint, \latin{e.g.} when there is few noisy data.
\item Because of this plug-in method, algorithm \ref{def.proc.gal} seems to be reasonable only if $\M_n$ is not too large, \latin{i.e.} if it has a polynomial complexity : $\card(\M_n) \leq \cM n^{\aM}$. Otherwise, we can for instance group the models of similar complexities and reduce $\M_n$ to a polynomial family.
\end{enumerate}
\end{remark}

\section{The histogram regression case} \label{sec.histos}
As studying algorithm~\ref{def.proc.gal} in general is a rather difficult question, we focus in this article on the case example of least-square regression on histograms. Although we do not consider histograms as a final goal, this first theoretical step will be useful to derive heuristics making the general algorithm~\ref{def.proc.gal} work.

We first precise the framework and some notations. The data $(X_i,Y_i) \in \X \times \R$ are i.i.d. of common law $P$. Denoting $\bayes$ the regression function, we have \begin{equation} \label{eq.donnees.reg} Y_i = \bayes(X_i) + \sigma(X_i) \epsilon_i \end{equation} where $\sigma : \X \mapsto \R$ is the heteroscedastic noise-level and $\epsilon_i$ are i.i.d. centered noise terms, possibly dependent from $X_i$, but with variance 1 conditionally to $X_i$.
The feature space $\X$ is typically a compact set of $\R^d$. We use the least-square contrast $\gamma : (t,(x,y)) \mapsto (t(x)-y)^2$ to measure the quality of a predictor $t : \X \mapsto \Y$. As a consequence, the Bayes predictor is the regression function $\bayes$, and the excess loss is $\perte{t} = \E_{(X,Y)\sim P} \carre{t(X)-\bayes(X)}$. 
To each model $S_m$, we associate the empirical risk minimizer $\ERM_m = \ERM_m(P_n) = \arg\min_{t \in S_m} \{ P_n \gamma(t) \}$ (when it exists and is unique).

Each model in $(S_m)_{\mM_n}$ is the set of piecewise constant functions (histograms) on some partition $(\Il)_{\lamm}$ of $\X$. It is thus a vector space of dimension $D_m = \card(\Lambda_m)$, spanned by the family $(\1_{\Il})_{\lamm}$. As this basis is orthogonal in $L^2(\mu)$ for any probability measure on $\X$, we can make explicit computations that will be useful to understand algorithm~\ref{def.proc.gal}. The following notations will be useful throughout this article.
\begin{equation*}
\pl \egaldef P(X\in \Il) \qquad \phl \egaldef P_n(X \in \Il) \qquad \phlW = \phl \Wl \egaldef \Pnb(X \in \Il) \end{equation*}
\begin{align*}
\bayes_m &\egaldef \arg\min_{t \in S_m} P \gamma(t) = \sum_{\lamm} \betl \1_{\Il} &\qquad \betl &= \E_P [ Y \sachant X \in \Il ] \\
\ERM_m &\egaldef \arg\min_{t \in S_m} P_n \gamma(t) = \sum_{\lamm} \bethl \1_{\Il} &\qquad \bethl &= \frac{1}{n\phl} \sum_{X_i \in \Il} Y_i \\
\ERMb_m &\egaldef \arg\min_{t \in S_m} \Pnb \gamma(t) = \sum_{\lamm} \bethlW \1_{\Il} &\qquad \bethlW &= \frac{1}{n\phlW} \sum_{X_i \in \Il} W_i Y_i
\end{align*}
Remark that $\ERM_m$ is uniquely defined if and only if each $\Il$ contains at least one of the $X_i$, and the same problem arises for $\ERMb_m$. This is why we will slightly modify the general algorithm for histograms. Before this, we compute the ideal penalty (assuming that $\min_{\lamm} \phl >0$ ; otherwise, the model $m$ should clearly not be chosen) :
\begin{equation*} 
\penid(m) = (P-P_n) \gamma(\ERM_m) = \sum_{\lamm} \paren{ \pl + \phl } \paren{  \bethl - \betl }^2 + (P-P_n) \gamma(\bayes_m) \enspace .
\end{equation*}
The last term in the sum being centered, it is estimated as zero by the resampling version of $\penid$. The first term is a sum of $D_m$ terms, each one depending only on the restrictions of $P$ and $P_n$ to $\Il$. Thus, if we assume that $\phl>0$ and if we compute separately all those terms, conditionally to $\phlW>0$, we can define the resampling version of $\penid(m)$. This leads to the following algorithm.
\begin{proc}[Resampling penalization for histograms] \label{def.proc.his}
\begin{enumerate}
\item[0.] Choose a threshold  $\seuilminphl \geq 1$ and replace $\M_n$ by 
\begin{equation*} 
\Mh_n = \left\{ \mM_n \telque \min_{\lamm} \{ n\phl \} \geq \seuilminphl \right\} \enspace .
\end{equation*}
\item[1.] Choose a resampling scheme $\mathcal{L}(W)$.
\item[2.] Choose a constant $C \geq C_W (\seuilminphl) $ where $C_W$ is defined by \eqref{def.CW}.
\item[3'.] Compute the following resampling penalty for each $\mMh_n$ :
\begin{equation*} 
\pen(m) = C \sum_{\lamm} \Es \croch{ \paren{ \phl + \phlW } \paren{  \bethlW - \bethl }^2 \Big\sachant \Wl >0} \enspace .
\end{equation*}
\item[4'.] Minimize the penalized empirical criterion to choose $\mh$ and thus $\ERM_{\mh}$ :
\begin{equation*} 
\mh \in \arg\min_{\mMh_n} \left\{ P_n \gamma(\ERM_m(P_n)) + \pen(m) \right\} \enspace .
\end{equation*} 
\end{enumerate}
\end{proc}
\begin{remark}
\begin{enumerate}
\item The two modifications of the algorithm for histograms do not affect much the result if $\seuilminphl$ is of the order $\ln(n)$. Indeed, models with very few data are not relevant in general, and if $\min_{\lamm} \{ n\phl\} \geq \seuilminphl$ is not too small, the event $\{ \Wl = 0 \}$ has a very small probability.
\item We allow $C$ to depend on $\seuilminphl$ since the ``optimal'' constant $C_W$ may depend on it, but this dependence is mild according to our computations.
\end{enumerate}
\end{remark}

When the resampling weights are exchangeable (see definition below), we are able to compute $\pen$ explicitly. It is enlightening to compare it with $\penid$ in expectation, conditionally to $(\phl)_{\lamm}$ (we denote by $\Em\croch{\cdot}$ this conditional expectation) :
\begin{align} \label{eq.Em_penid}
\Em\croch{ \penid(m) } &= \frac{1}{n} \sum_{\lamm} \paren{ 1 + \frac{\pl}{\phl} } \paren{\carre{\sigla} + \carre{\sigld}} \\ \notag
\Em\croch{ \pen(m) } &= \frac{C}{n} \sum_{\lamm} \paren{ R_{1,W}(n,\phl) + R_{2,W}(n,\phl) } \paren{\carre{\sigla} + \carre{\sigld}} \end{align}
\begin{equation*}
\mbox{with} \quad  \carre{\sigla} \egaldef \E[\sigma(x)^2 \sachant X \in \Il] \enspace \mbox{;} \quad \carre{\sigld} \egaldef \E[(\bayes(X)-\bayes_m(X))^2 \sachant X \in \Il] \end{equation*} 
\begin{align*}
\mbox{and for $k=1,2$} \qquad R_{k,W}(n,\phl) &= \E\croch{ \frac{(W_i - \Wl)^2}{\Wl^{3-k}} \Big\sachant \Wl>0} \enspace .
\end{align*}
Hence, contrary to Mallows' penalty (with $\sigma^2$ known or estimated), resampling penalties really take into account the heteroscedasticity of the noise ($\sigla$ depends on $\lambda$) and the bias terms $\carre{\sigld}$. We then define
\begin{equation} \label{def.CW}
C_W(\seuilminphl) \egaldef \sup_{n\phl \geq \seuilminphl} \left\{ \frac{2}{R_{1,W}(n,\phl) + R_{2,W} (n,\phl)}\right\} 
\end{equation}
and $C_W^{\prime}(\seuilminphl)$ is the infimum of the same quantity.

\subsection*{Examples of resampling weights}
In this article, we consider resampling weights $W = (W_1,\ldots, W_n) \in [0;+\infty)^n$ such that $\E[W_i]=1$ for all $i$ and $\E[W_i^2]<\infty$. We mainly consider the following exchangeable weights (\latin{i.e.} such that for any permutation $\tau$, $(W_{\tau(1)}, \ldots, W_{\tau(n)}) \egalloi (W_{1}, \ldots, W_{n})$).
\begin{enumerate}
\item {\em Efron} ($q$) : multinomial vector with parameters $(q;n^{-1},\ldots, n^{-1})$. Then, $R_{2,W}(n,\phl) = (n/q)\times(1-(n\phl)^{-1})$. A classical choice is $q=n$.
\item {\em Rademacher} : $W_i$ i.i.d., 2 times Bernoulli($1/2$). Then, $R_{2,W}(n,\phl) = 1$.
\item {\em Random hold-out} ($q$) (or cross-validation) : $W_i = \frac{n}{q} \1_{i \in I}$ with $I$ uniform random subset (of cardinality $q$) of $\{1, \ldots, n\}$. $R_{2,W}(n,\phl)=(n/q)-1$. A classical choice is $q=n/2$.
\item {\em Leave-one-out} = Random hold-out ($n-1$). Then, $R_{2,W}(n,\phl) = (n-1)^{-1}$.
\end{enumerate}
In each case, we can show that $R_{1,W} = R_{2,W}(1 + \delta^{(W)}_{n,\phl})$ for some explicit small term $\delta^{(W)}_{n,\phl}$ (numerically of the same order as $\E[\pl/\phl\sachant \phl>0]-1$ in expectation for the three first resamplings, and slightly smaller in the Leave-one-out case). Thus, $C_W \approx C_W^{\prime} \approx R_{2,W}^{-1}$ (asymptotically in $\seuilminphl$).
 
For computational reasons, it is also convenient to introduce the following {\em $V$-fold cross-validation} resampling weights : given a partition $(B_j)_{1 \leq j \leq V}$ of $\{1, \ldots, n\}$ and $W^B \in \R^V$ leave-one-out weights, we define $W_i = W^B_j$ for each $i \in B_j$. The partition should be taken as regular as possible, and then we can compute $\E[\pen(m)]$ and show that $C_W \approx V-1$.

The Rademacher weights lead to penalties close in spirit to local Rademacher complexities (the link between global Rademacher complexities and global resampling penalties with Rademacher weights can be found in \cite{Fro:2004}). The links with the classical leave-one-out and VFCV algorithms are given in Sect.~\ref{sec.V-fold}. 

\section{Main results} \label{sec.main}
In this section, we prove that algorithm~\ref{def.proc.his} has some optimality properties under the following restrictions for some non-negative constants $\aM$, $\cM$, $c_A$, $c_{\mathrm{rich}}$~:
\begin{enumerate}
\item[\hyp{P1}] \label{hyp.proc.Mpoly} Polynomial complexity of $\M_n$ : $\card(\M_n) \leq \cM n^{\aM}$.
\item[\hyp{P2}] \label{hyp.proc.Mrich} Richness of $\M_n$ : $\forall x \in [1,n c_{\mathrm{rich}}^{-1}]$, $\exists \mM_n$ s.t. $D_m \in [x;c_{\mathrm{rich}}x]$.
\item[\hyp{P3}] \label{hyp.proc.exch} The weights are exchangeable, among the examples given in Sect. \ref{sec.histos}. 
\item[\hyp{P4}] \label{hyp.proc.An} The threshold is large enough : $C_{A} \ln(n) \geq \seuilminphl \geq (26 + 7\aM) \ln(n)$.
\end{enumerate}
Assumption \hyp{P1} is almost necessary, since too large families of models need larger penalties than polynomial families \cite{Bir_Mas:2002,Bar:2002,Sau:2006}.
Assumption \hyp{P2} is necessary but it is always satisfied in practice.
Assumption \hyp{P3} is only here to ensure that we have an explicit formula for the penalty, and sharp bounds on $R_{1,W}$ and $R_{2,W}$.
The constant $(26 + 7\aM)$ in \hyp{P4} is quite large due to technical reasons, but much smaller values (larger than 2) should suffice in practice.

\begin{theorem} \label{the.oracle_traj_non-as}
Assume that the $(X_i,Y_i)$'s satisfy the following assumptions :
\begin{enumerate}
\item[\hyp{Ab}] Bounded data : $\norm{Y_i}_{\infty} \leq A < \infty$.
\item[\hyp{An}] Noise-level bounded from below : $\sigma(X_i) \geq \sigmin>0$ a.s.
\item[\hyp{Ap}] Polynomial decreasing of the bias : \[ \exists \beta_1\geq\beta_2>0, \, C_{\bayes},c_{\bayes}>0 \quad \mbox{s.t.} \quad c_s D_m^{-\beta_1} \leq \perte{\bayes_m} \leq C_s D_m^{-\beta_2} \enspace . \]
\item[\hyp{Ar}] (pseudo)-Regular histograms : $\forall \mM_n$, $\min_{\lamm} \{\pl\} \geq c_{\mathrm{reg}} D_m^{-1}$. 
\end{enumerate}

Let $\mh$ be the model chosen by algorithm~\ref{def.proc.his} (under restrictions \hyp{P1--4}), with $\eta^{\prime} C_W^{\prime}(\seuilminphl) \geq C \geq \eta C_W(\seuilminphl)$ for some $\eta,\eta^{\prime}>\frac{1}{2}$. It satisfies, with probability at least $1 - L_{\mhyp{A},\mhyp{P}} n^{-2}$ ($L_{\mhyp{A},\mhyp{P}}$ may depend on constants in \hyp{A} and \hyp{P}, but not on $n$), 
\begin{equation} \label{eq.oracle_traj_non-as}
\perte{\ERM_{\mh}} \leq K(\eta,\eta^{\prime}) \inf_{\mM_n} \left\{ \perte{\ERM_m} \right\} \enspace .
\end{equation}
At the price of enlarging $L_{\mhyp{A},\mhyp{P}}$, the constant $K(\eta,\eta^{\prime})$ can be taken close to $(1+2(\eta^{\prime}-1)_+) (1-2(1-\eta)_+)^{-1}$, where $x_+ \egaldef \max(x,0)$. In particular, $K(\eta,\eta^{\prime})$ is almost 1 if $\eta$ and $\eta^{\prime}$ are close to 1.

Moreover, we have the oracle inequality
\begin{equation} \label{eq.oracle_class_non-as}
\E \croch{\perte{\ERM_{\mh}}} \leq K(\eta,\eta^{\prime}) \E \croch{ \inf_{\mM_n} \left\{ \perte{\ERM_m} \right\} } + \frac{A^2 L_{\mhyp{A},\mhyp{P}}}{n^2} \enspace .
\end{equation}
\end{theorem}
\begin{proof}[sketch]
By definition of $\mh$,
\[ \forall \mMh_n, \quad (\pen - \penid^{\prime})(\mh) + \perte{\ERM_{\mh}} \leq \perte{\ERM_m} + (\pen - \penid^{\prime})(m) \]
where we replaced $\penid$ by $\penid^{\prime} \egaldef \penid - (P_n - P)\gamma(\bayes)$. 
In order to obtain \eqref{eq.oracle_traj_non-as} with $\Mh_n$ instead of $\M_n$, we show concentration inequalities for $\pen(m)-\penid^{\prime}(m)$ around zero, with remainders $\ll \perte{\ERM_m}$ if $D_m$ is large (larger than some power of $\ln(n)$). We use the following steps :
\begin{enumerate}
\item explicit computation of $\penid^{\prime}$ and $\pen$ when $W$ is exchangeable.
\item accurate bounds on $R_{1,W}$ and $R_{2,W}$, so that $(1-\delta(\seuilminphl)) \Em[2 p_2(m)] \leq \Em[\pen(m)] \leq (1+\delta(\seuilminphl)) \Em[2 p_2(m)]$ with $p_2(m) = P_n \paren{\gamma(\bayes_m) - \gamma(\ERM_m)}$ and $\lim_{\seuilminphl \rightarrow \infty} \delta = 0$. This needs sharp bounds on $\E[Z^{-1} \sachant Z>0]$ with $\mathcal{L}(Z) = \mathcal{L}(\Wl \sachant \phl)$, for each resampling scheme introduced in Sect.~\ref{sec.histos}.
\item moment inequalities for $\pen$, $p_2$ and $p_1(m) = P \paren{\gamma(\ERM_m) - \gamma(\bayes_m)}$, conditionally to $(\phl)_{\lamm}$, around their conditional expectations. This step uses results from \cite{Bou_Bou_Lug_Mas:2005}, or can be derived from \cite{Gin_Lat_Zin:2000}, since all those quantites are U-statistics of order 2 (this last fact is not true without the conditioning). This implies (unconditional) concentration inequalities.
\item concentration inequality for $(P_n - P)(\gamma(\bayes_m) - \gamma(\bayes))$ (Bernstein's inequality suffices in the bounded case).
\item since $\Em \croch{p_2(m)} = \E \croch{p_2(m)}$, it only remains to prove that $\Em \croch{p_1(m)} \approx \E \croch{p_1(m)}$ and $p_2 \approx p_1$ with high probability. We here use the Cram\'er-Chernoff method (it can be used since the $(\phl)_{\lamm}$ are negatively associated \cite{Dub_Ran:1998}), together with estimates of the exponential moments of the inverse of a binomial random variable. Controlling the remainder needs a lower bound on $\min_{\lamm} \{n\pl\}$ that comes from \hyp{P4} (and Bernstein's inequality).
\item using the assumptions, all the remainders in our concentration inequalities are much smaller than $\E[\perte{\ERM_{\mh}}]$ when $D_m \geq D_0(n) = c_1 (\ln(n))^{c_2}$ (with $c_1,c_2$ depending on the constants in the assumptions). 
\end{enumerate}

Let $m^{\star}$ be a minimizer of $\perte{\ERM_{m}}$ over $\M_n$ (with an infinite loss when $\ERM_m$ is not uniquely defined). It remains to prove that, with large probability, $D_{\mh} \geq D_0(n)$, $D_{m^{\star}} \geq D_0(n)$ and $m^{\star} \in \Mh_n$.
These hold for $n$ large enough thanks to \hyp{Ap} and \hyp{Ar} (we did not use \hyp{Ar} before).

We finally show that \eqref{eq.oracle_traj_non-as} implies \eqref{eq.oracle_class_non-as} : let $\Omega_n$ be the event of probability $1-L_{\mhyp{A},\mhyp{P}}n^{-2}$ on which \eqref{eq.oracle_traj_non-as} occurs. On $\Omega_n^c$, $\perte{\ERM_{\mh}}$ is bounded by $A^2$, so that
\begin{align*} \E \croch{\perte{\ERM_{\mh}}} &=  \E \croch{\perte{\ERM_{\mh}} \1_{\Omega_n}} + \E \croch{\perte{\ERM_{\mh}} \1_{\Omega_n^c}} \\
&\leq K(\eta,\eta^{\prime}) \E \croch{\inf_{\mM_n} \perte{\ERM_m}} + L_{\mhyp{A},\mhyp{P}} A^2 n^{-2} \enspace .  \quad \qed \end{align*} 
\end{proof}

Theorem \ref{the.oracle_traj_non-as} implies the a.s. asymptotic optimality of algorithm~\ref{def.proc.gal.pente} in this framework. This means that if $\bayes$ and $\sigma(X)$ do not make the model selection problem too hard, the resampling penalization algorithm is working, without any knowledge on the smoothness of $\bayes$, the heteroscedasticity of $\sigma$ or any property that the unknown law $P$ may satisfy. In that sense, it is a {\em naturally adaptive algorithm}.

\medskip

The lower bound in assumption \hyp{Ap} may seem strange, but it is intuitive that when the bias is decreasing very fast, the optimal model is of quite small dimension. Then, bounds relying on the fact that this dimension is large can not work. The same kind of assumption has already been used in the density estimation framework for the same reason \cite{Sto:1985}. 

Moreover, we can prove that non-constant h\"olderian functions satisfy \hyp{Ap} when $X$ has a lower-bounded density w.r.t. the Lebesgue measure on $\X \subset \R$. The following result states that resampling penalization is adaptive to the h\"olderian smoothness of $\bayes$ in an heteroscedastic framework, since it attains the minimax rate of convergence $n^{-2\alpha/(2\alpha+1)}$ \cite{Sto:1980}.
\begin{theorem} \label{the.holder} Let $\X$ be a compact interval of $\R$ and $\Y\subset \R$. Assume that $(X_i,Y_i)$ satisfy \hyp{Ab}, \hyp{An} and the following assumptions :
\begin{enumerate}
\item[\hyp{Ad}] Density bounded from below : $\exists c_{\min}^X>0$, $\forall I \subset \X$, $P(X \in I) \geq c_{\min}^X \Leb(I)$.
\item[\hyp{Ah}] H\"olderian regression function : there exists $\alpha\in(0;1]$ and $R>0$ s.t. 
\begin{equation*} \bayes \in \mathcal{H}(\alpha,R) \quad \mbox{ \latin{i.e.} } \quad \forall x_1, x_2 \in \X, \, \absj{\bayes(x_1) - \bayes(x_2)} \leq R \absj{x_1 - x_2}^{\alpha} \enspace . \end{equation*}
\end{enumerate} 

Let $\M_n$ be the family of regular histograms of dimensions $1 \leq D \leq n$, $\mh$ the model chosen by algorithm~\ref{def.proc.his}, with \hyp{P3-4} satisfied ($\aM=0$) and $C$ like in Theorem~\ref{the.oracle_traj_non-as}.
Then, denoting $\sigmax = \sup_{\X} \absj{\sigma} \leq A$, there are some constants $L_{2,\mhyp{A},\mhyp{P}}$ (that may depend on all the constants in the assumptions) and $L_1(\eta,\eta^{\prime},\alpha)$ such that
\begin{equation} \label{eq.oracle_class_non-as.hold}
\E \croch{\perte{\ERM_{\mh}}} \leq L_1 n^{-2\alpha/(2\alpha+1)} R^{2\alpha/(2\alpha+1)} \sigmax^{4\alpha/(2\alpha+1)} + L_{2,\mhyp{A},\mhyp{P}} n^{-2} \enspace .
\end{equation}
Moreover, if $\sigma$ is $K_{\sigma}$-Lipschitz, the constant $\sigmax^2$ may be replaced by $\int_{\X} \sigma(t)^2 dt$ (at the price of enlarging $L_{2,\mhyp{A},\mhyp{P}}$).
\end{theorem} 
\begin{proof}[sketch]
\begin{enumerate}
\item Since $\alpha \in (0;1]$, any non-constant function $\bayes \in \mathcal{H}(\alpha,R)$ satisfies \hyp{Ap} with $\beta_2 = 2\alpha$ and $\beta_1=1+\alpha^{-1}$ (the lower bound uses \hyp{Ad}).
\item Assumptions \hyp{P1}, \hyp{P2} and \hyp{Ar} are automatically satisfied by the regular family, so we can use \eqref{eq.oracle_class_non-as}. From the proof of Theorem~\ref{the.oracle_traj_non-as}, we obtain estimations of $\E\croch{\perte{\ERM_m}}$. Optimizing in $D_m$ gives \eqref{eq.oracle_class_non-as.hold} for non-constant functions.
\item When $\bayes$ is constant, a direct proof shows that $D_{\mh}$ is at most of order $\ln(n)^{\xi_1}$ with large probability. This ensures that $\E\croch{\perte{\ERM_{\mh}}}$ is at most of order $(\ln(n))^{\xi_2} n^{-1} \ll n^{-2\alpha / (2\alpha +1)}$ for every $\alpha>0$. \qed
\end{enumerate}
\end{proof}

Other results like Theorem~\ref{the.oracle_traj_non-as} may be proved under other assumptions : unbounded data (with moment inequalities for the noise, regularity assumptions on $\bayes$ and an upper bound on $\sigma$), $\sigma(x)$ that can vanish (with the unbounded assumptions, $\E[\sigma^2(X)]>0$ and some regularity on $\sigma$), etc. We skip their detailed statements in order to focus on the last two sections, where we give a new look on $V$-fold cross-validation (seen from the penalization viewpoint) and illustrate theoretical results with a simulation study.

\section{Links with $V$-fold cross-validation} \label{sec.V-fold}
The results of Sect.~\ref{sec.main} assume that the weights are exchangeable. However, computing exactly the resampling penalties with such weights may be quite long~: without a closed formula for $\pen$, $\ERMb_m$ has to be computed for at least $n$ (and up to $2^n$) different weight vectors. Using the $V$-fold idea, we defined VFCV weights in Sect.~\ref{sec.histos}, that allows to compute each penalty by considering only $V$ different weight vectors. We call the resulting algorithm penVFCV.

It is quite enlightening to compare penVFCV to a more classical version of VFCV, where the final estimator is $\ERM_{\mh}$ with 
\begin{equation} \label{def.mh.VFCvclass}
\mh \in \arg\min_{\mM_n} \left\{  \crit_{\mathrm{VFCV}}(m) \right\} = \arg\min_{\mM_n} \left\{  \frac{1}{V} \sum_{j=1}^V P_n^{(j)} \gamma(\ERM_m^{(-j)}) \right\} \enspace .
\end{equation}
The superscript $(j)$ (resp. $(-j)$) above means that $P_n$ and $\ERM_m$ are computed with the data belonging to the block $B_j$ (resp. to $B_j^c$).
Assuming that the $V$ blocks have the same size (and forgetting unicity issues of $\ERM_m^{(-j)}$, that may be solved as before), we have (for any $j$)
\begin{align} \label{eq.Ecrit.VFCVclass}
\E \croch{ \crit_{\mathrm{VFCV}}(m) } &= P \gamma(\bayes_m) + \E \croch{ P \gamma(\ERM_m^{(-j)}) - P \gamma(\bayes_m)} \\ \notag
&= P \gamma(\bayes_m) + \frac{V}{(V-1)n} \sum_{\lamm} \paren{1 + \delta^{(V)}_{n,\pl}} \paren{ \carre{\sigla} + \carre{\sigld} } \enspace \end{align}
where $\delta^{(V)}_{n,\pl}$ is typically small and non-negative (when $n\pl$ is large enough).

On the other hand, 
we can compute exactly the expectation of the penVFCV criterion (with a constant $C=C_W=V-1$) when the blocks have the same size~:
\begin{equation} \label{eq.Ecrit.VFCVpen}
\E \croch{ \crit_{\mathrm{penVFCV}}(m) } = P \gamma(\bayes_m) + \frac{1}{n} \sum_{\lamm} \paren{1 + \delta^{(penV)}_{n,\pl}} \paren{ \carre{\sigla} + \carre{\sigld} }  \enspace \end{equation}
for some typically small non-negative $\delta^{(penV)}_{n,\pl}$. 

\medskip

Comparing \eqref{eq.Ecrit.VFCVclass} and \eqref{eq.Ecrit.VFCVpen} with \eqref{eq.Em_penid}, one can see that up to small terms, both criterions are in expectation the sum of the bias and a variance term. The main difference between them lies in the constant in front of the variance : it is equal to $C/(V-1)=1$ for penVFCV, whereas it is equal to $V/(V-1)>1$ for VFCV. 

The classical $V$-fold cross-validation is thus ``overpenalizing'' within a factor $V/(V-1)$ because it estimates the generalization ability of $\ERM^{(-j)}_m$, which is built upon less data than $\ERM_m$. This enlightens some clues for the choice of $V$ : {\em computational issues} (the smaller $V$, the faster will be the algorithm), {\em stability} of the algorithm ($V=2$ is known to be quite unstable, and leave-one-out much more stable), and {\em overpenalization} ($V/(V-1)$ should not be too far from 1). Our analysis do not quantify the stability issue, but it is sufficient to explain why the asymptotic optimality of leave-$p$-out needs $p \ll n$ for a prediction purpose \cite{KCLi:1987} and $p \sim n$ for an identification purpose \cite{Zha:1993}. Indeed, the overpenalization factor is $n/(n-p) = (1-p/n)^{-1}$  should go to 1 for optimal prediction and to infinity for a.s. identification. Moreover, from the non-asymptotic viewpoint ($n$ small and $\sigma$ large, or $\bayes$ irregular), it is known that overpenalization (\latin{i.e.} positively biased penalties) gives better results. This means that the better $V$ may not always be the largest one for classical $V$-fold, independently from computational issues.

On the contrary, penVFCV is not overpenalizing, unless we explicitly choose $C > C_W$. We thus do not have to take into account the third factor for choosing $V$, so that it may be more accurate than VFCV within a smaller computation time. In the non-asymptotic viewpoint (or for an identification purpose), it is also easier to overpenalize when we need to, without destabilizing the algorithm by taking a small $V$.

\medskip

A refined analysis of the ``negligible'' terms such as $\delta^{(penV)}_{n,\pl}$, compared to the expectation of $\pl/\phl$, explains why the leave-one-out may be overfitting a little (see the simulations hereafter). We do not detail this phenomenon since it disappears when $V/(V-1)$ stays away from 1.

\section{Simulations} \label{sec.simus}
To illustrate the results of Sect. \ref{sec.main} and the analysis of Sect. \ref{sec.V-fold}, we compare the performances of algorithm~\ref{def.proc.his} (with several resampling schemes),  Mallows' $C_p$ and VFCV on some simulated data. 

\begin{figure}
\begin{minipage}[b]{.32\linewidth}   \centerline{\epsfig{file=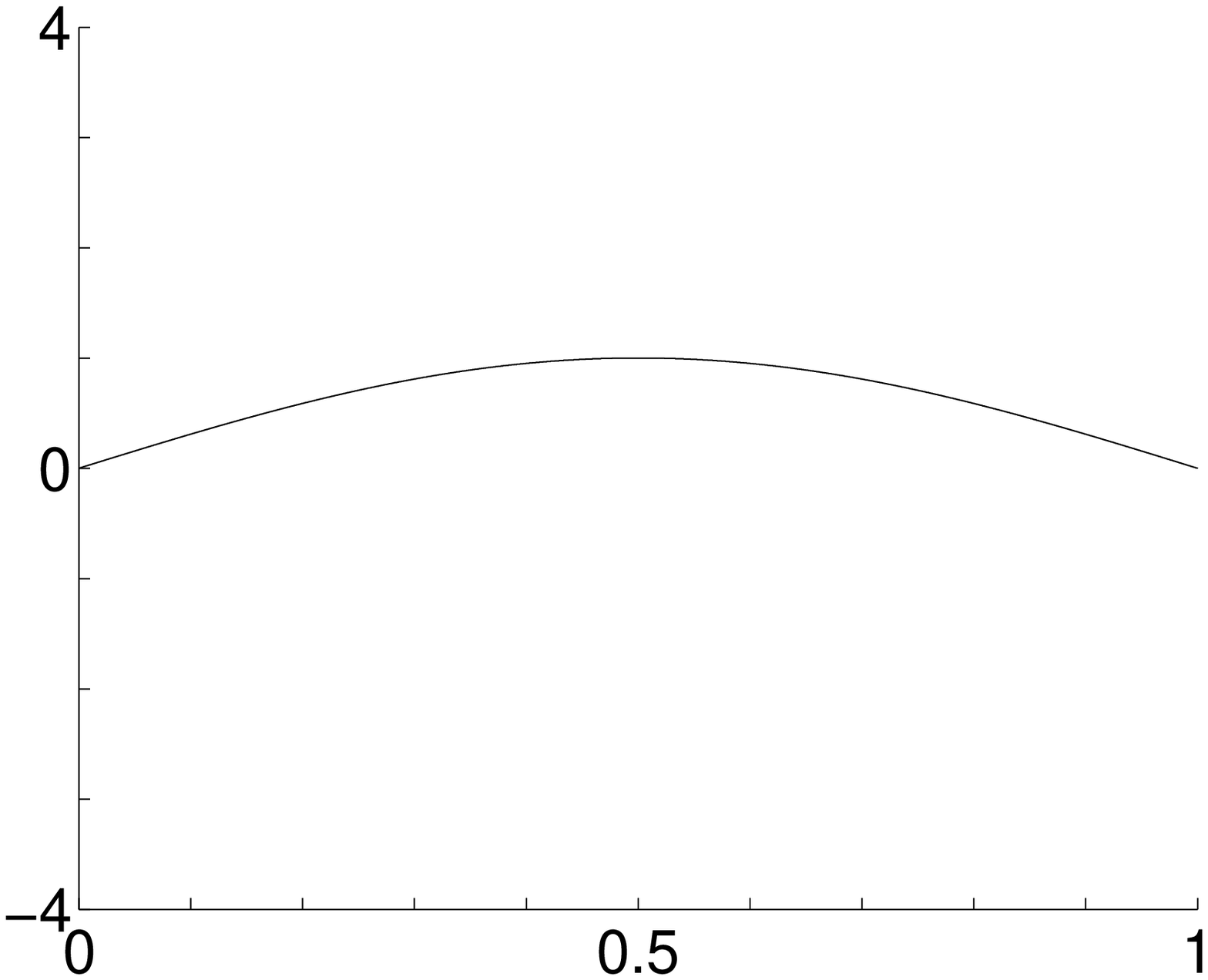,width=0.9\textwidth}}
   \caption{$\bayes(x) = \sin(\pi x)$\label{fig.sin.fonc}}
\end{minipage} \hfill
\begin{minipage}[b]{.32\linewidth}   \centerline{\epsfig{file=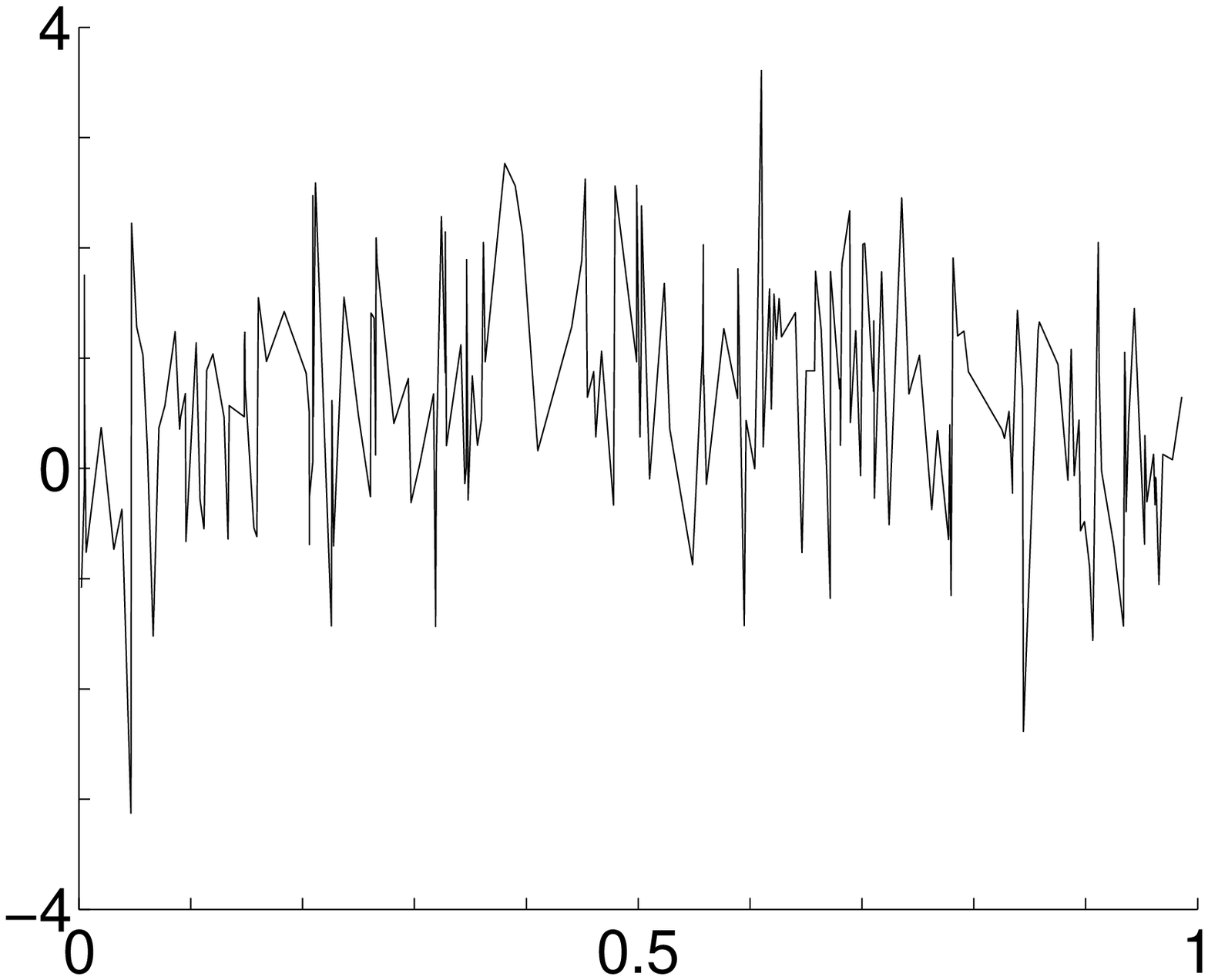,width=0.9\textwidth}}
   \caption{S1}
   \label{fig.S1.data}
\end{minipage} \hfill
 \begin{minipage}[b]{.32\linewidth}   \centerline{\epsfig{file=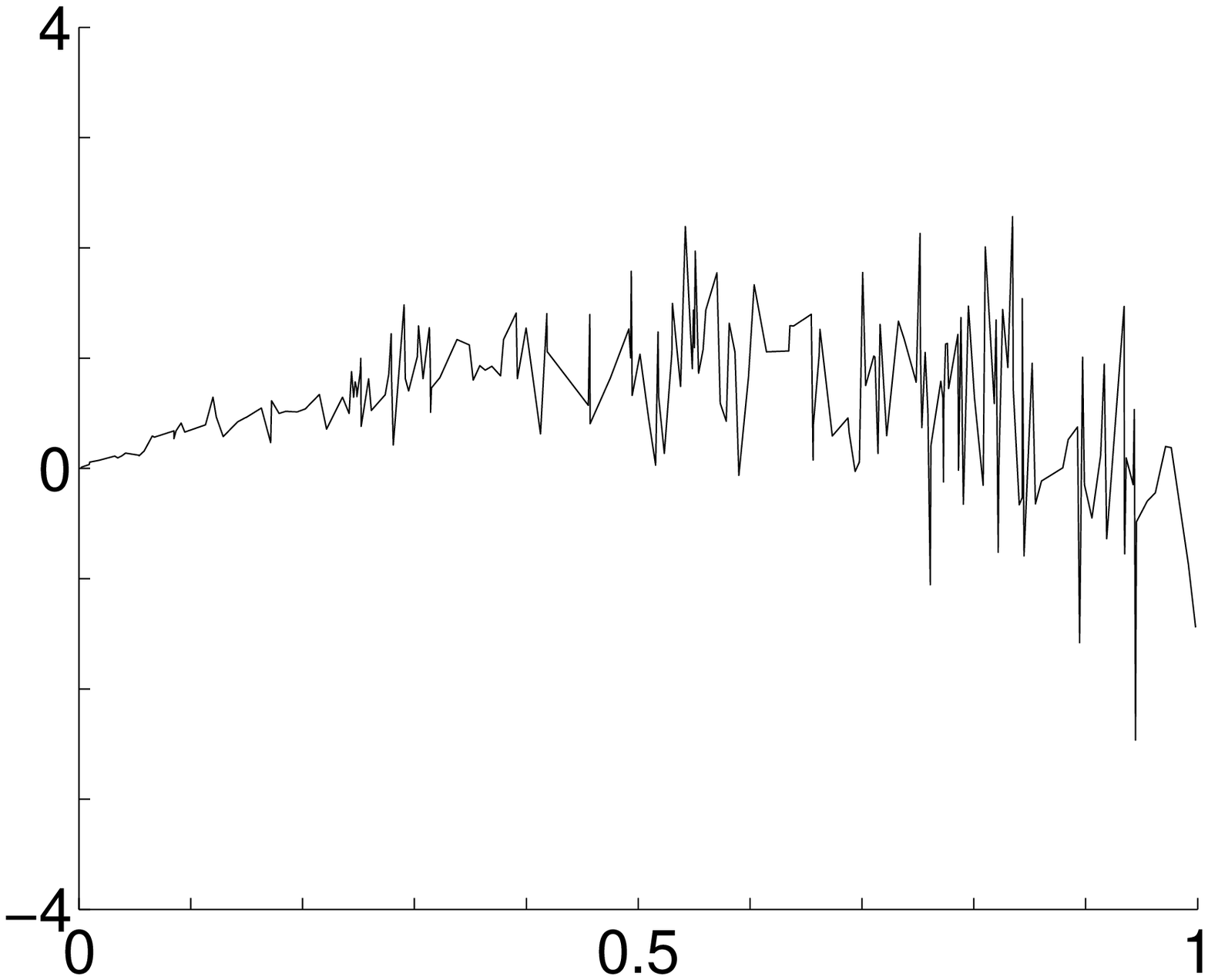,width=0.9\textwidth}}
   \caption{S2}
   \label{fig.S2.data}
\end{minipage}
\end{figure}

\begin{figure}
 \begin{minipage}[b]{.32\linewidth}   \centerline{\epsfig{file=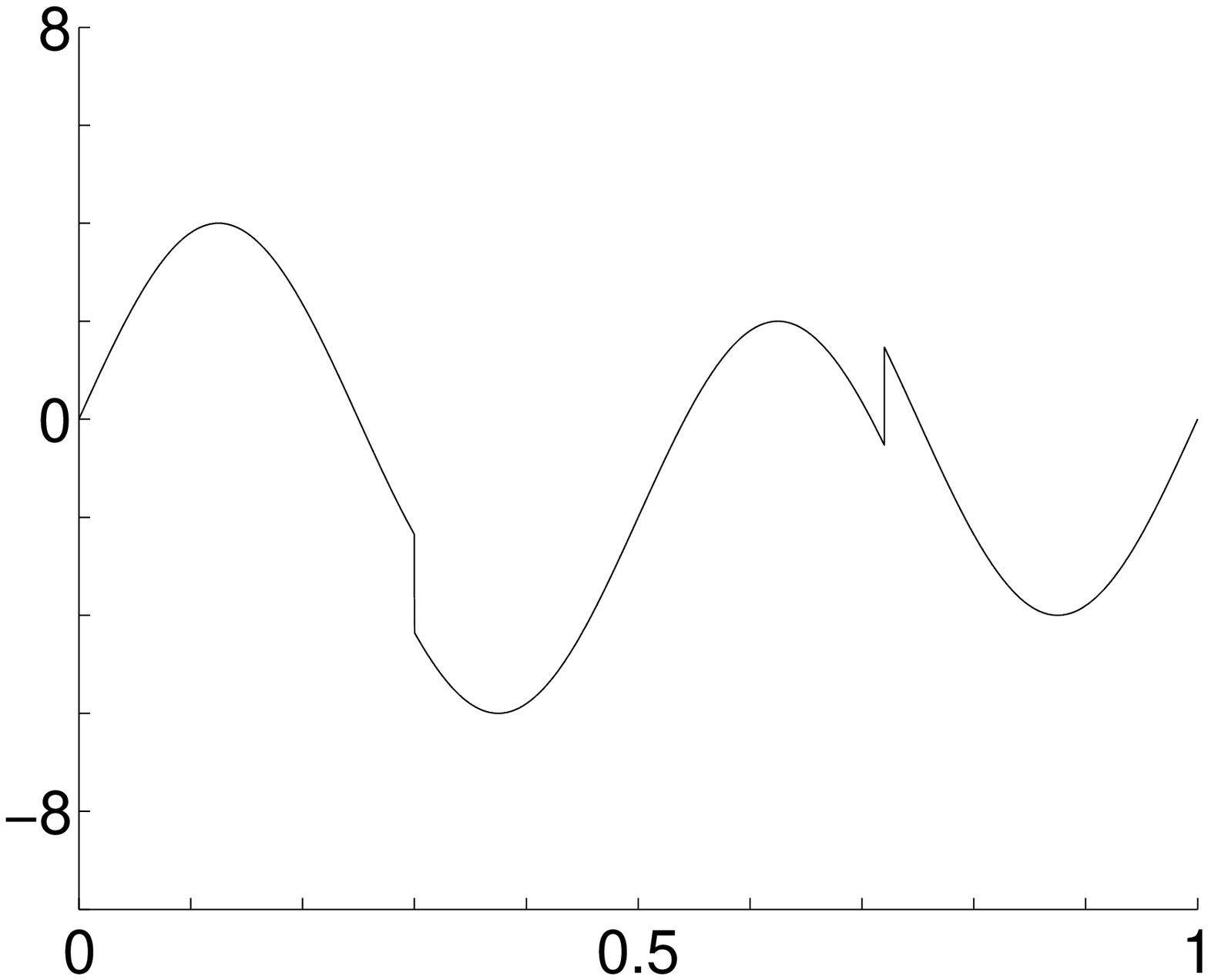,width=0.9\textwidth}}
   \caption{$\bayes = \mathrm{HeaviSine}$ (see \cite{Don_Joh:1995})\label{fig.Hsin.fonc}}
\end{minipage}
\begin{minipage}[b]{.32\linewidth}   \centerline{\epsfig{file=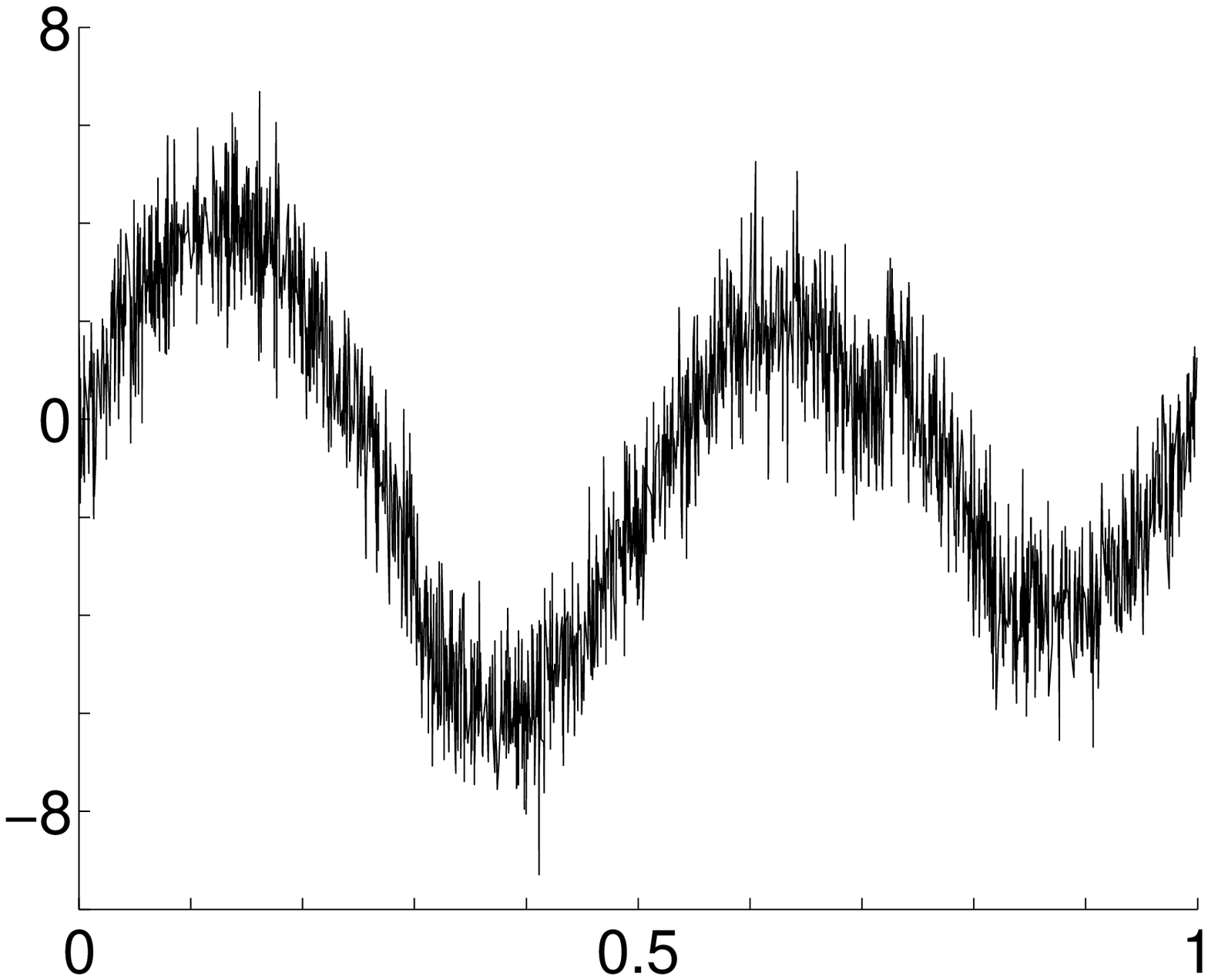,width=0.9\textwidth}}
   \caption{HSd1}
   \label{fig.HSd1.data}
\end{minipage} \hfill
 \begin{minipage}[b]{.32\linewidth}   \centerline{\epsfig{file=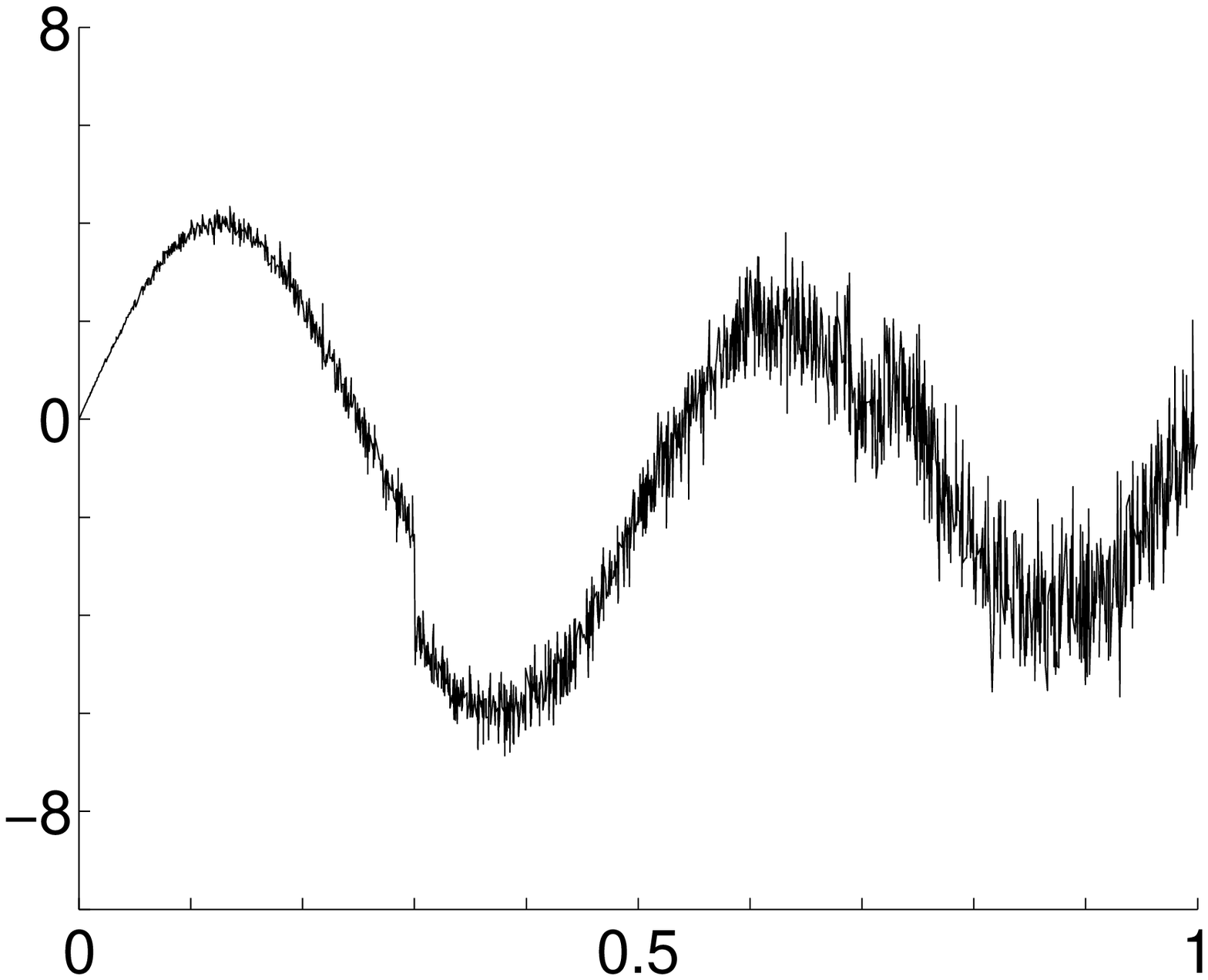,width=0.9\textwidth}}
   \caption{HSd2}
   \label{fig.HSd2.data}
\end{minipage}
\end{figure}

We report here four experiments, called S1, S2, HSd1 and HSd2. Data are generated according to \eqref{eq.donnees.reg} with $X_i$ i.i.d. uniform on $\X=[0;1]$ and $\epsilon_i \sim \mathcal{N}(0,1)$ independent from $X_i$.
They differ from the regression function $\bayes$ (smooth for S, see Fig.~\ref{fig.sin.fonc} ; smooth with jumps for HS, see Fig.~\ref{fig.Hsin.fonc}), the noise type (homoscedastic for S1 and HSd1, heteroscedastic for S2 and HSd2), the number $n$ of data, and are repeated $N = 1000$ times. Instances of data sets are given in Fig.~\ref{fig.S1.data}-\ref{fig.S2.data} and~\ref{fig.HSd1.data}-\ref{fig.HSd2.data}.
Their last difference lies in the families of models $\M_n$ :
\begin{enumerate}
\item[S1] regular histograms with $1 \leq D \leq \frac{n}{\ln(n)}$ pieces.
\item[S2] histograms regular on $\croch{0;\frac{1}{2}}$ and on $\croch{\frac{1}{2};1}$, with $D_1$ (resp. $D_2$) pieces, $1 \leq D_1,D_2 \leq \frac{n}{2\ln(n)}$. The model of constant functions is added to $\M_n$.
\item[HSd1] dyadic regular histograms with $2^k$ pieces, $0 \leq k \leq \ln_2(n) - 1$.
\item[HSd2] dyadic regular histograms with bin sizes $2^{-k_1}$ and $2^{-k_2}$, $0 \leq k_1,k_2 \leq \ln_2(n) - 1$ (dyadic version of S2). The model of constant functions is added to $\M_n$.
\end{enumerate}

We compare the following algorithms : 
\begin{enumerate}
\item[Mal] Mallows' $C_p$ penalty : $\pen(m) = 2 \widehat{\sigma}^2 D_m n^{-1}$ where $\widehat{\sigma}^2$ is the variance estimator used in \cite{Bar:2000}, Sect.~6.
\item[VFCV] Classical $V$-fold cross-validation, defined by \eqref{def.mh.VFCvclass}, with $V \in \{ 2,5,10,20 \}$.
\item[penEfr] Efron ($n$) penalty, $C= C_W=1$.
\item[penRad] Rademacher penalty, $C= C_W=1$.
\item[penRHO] Random hold-out ($n/2$) penalty, $C= C_W=1$.
\item[penLOO] Leave-one-out penalty, $C= C_W=n-1$.
\item[penVFCV] $V$-fold penalty, with $V \in \{ 2,5,10,20 \}$. $C= C_W=V-1$.
\end{enumerate}
For each of these except VFCV, we also consider the same penalties multiplied by $5/4$ (denoted by a $+$ symbol added after its shortened name). This intends to test for overpenalization.

In each experiment, for each simulated data set, we first remove the models with less than $\seuilminphl=2$ data points in one piece of their associated partition.  Then, we compute the least-square estimators $\ERM_m$ for each $\mMh_n$. Finally, we select $\mh \in \Mh_n$ using each algorithm and compute its true excess risk $\perte{\ERM_{\mh}}$ (and the excess risk of each model $m \in \M_n$). 
Since we simulate $N$ data sets, we can then estimate the two following benchmarks~:
\begin{equation*} 
C_{\mathrm{or}} = \frac{ \E\croch{ \perte{\ERM_{\mh}} }} {\E\croch{ \inf_{\mM_n}  \perte{\ERM_m} }} \qquad C_{\mathrm{path-or}} = \E\croch{\frac{  \perte{\ERM_{\mh} }} {\inf_{\mM_n} \perte{\ERM_m} } } 
\end{equation*}
Basically, $C_{\mathrm{or}}$ is the constant that should appear in an oracle inequality like \eqref{eq.oracle_traj_non-as}, and $C_{\mathrm{path-or}}$ corresponds to a pathwise oracle inequality like \eqref{eq.oracle_class_non-as}. As $C_{\mathrm{or}}$ and $C_{\mathrm{path-or}}$ approximatively give the same rankings between algorithms, we only report $C_{\mathrm{or}}$ in Tab.~\ref{tab.tout}. 

\medskip

We always observe that penRad and penRHO are competitive with Mal (S1) and much better for more ``difficult'' problems (S2 is heteroscedastic ; jumps in HSd1 and HSd2 induce much bias). On the other hand, VFCV is a little worse than Mal for easy problems (S1) and better for more difficult ones, but never better than penRad or penRHO.

The best resampling schemes (not taking overpenalization into account) are penRad and penRHO, in view of S1 and S2 (dyadic models do not induce much differences between them in HSd1 and HSd2). Then, penLOO is slightly underpenalizing and penEfr strongly overfits. The comparison penRad $\approx$ penRHO $>$ penLOO $\gg$ penEfr can also be derived from Sect.~\ref{sec.histos}. 

In the four experiments, overpenalizing within a factor $5/4$ leads to better results, 
mainly because $n$ is quite small for the noisy (S1, S2) or irregular (HSd1, HSd2) signals observed. This is no longer the case for some larger $n$ or smaller~$\sigma$. 

\medskip

We consider now $V$-fold algorithms. VFCV is slightly better than penVFCV, but worse than penVFCV+. The influence of $V$ on $C_{\mathrm{or}}$ confirms the discussion of Sect.~\ref{sec.V-fold}. For VFCV, the best $V$ may be $V=2$ (which overpenalizes, HSd1) or $V=20$ (which is more stable, HSd2), or even both (S1,S2). On the contrary, penVFCV (and penVFCV+) is always improved when $V$ increases, or at least it does not get worse. Then, the best one is penLOO (or penLOO+), \latin{i.e.} $V=n$, the small terms $\delta_{n,\pl}^{(penV)}$ being far less important than stability. 
This enlightens the interest of defining $V$-fold penalties, for which it is easier to solve the complexity-accuracy trade-off.

\begin{remark}
We only report here the result of 4 experiments, but several other ones (with $n$ larger, $\sigma$ smaller, $\sigma(x)=\1_{x \in \croch{\frac{1}{2};1}}$ or other regression functions $\bayes$ such as Doppler, $\sqrt{\cdot}$ and a regular histogram) give the same kind of results. The constants $C_{\mathrm{or}}$ and $C_{\mathrm{path-or}}$ are decreasing to 1 when $n$ increases and $\sigma$ decreases. 

The overpenalization factor $5/4$ is generally not optimal, and even not always better than $1$ (in particular when $n$ is large or $\sigma$ small). We have for instance $C_{\mathrm{or}}(\mathrm{penLOO}) < C_{\mathrm{or}}(\mathrm{penRHO}) < C_{\mathrm{or}}(\mathrm{penRHO+})$ in S1 with $\sigma \equiv 0.1$ (with only small differences).
\end{remark}

\begin{table} 
\caption{Accuracy indexes $C_{\mathrm{or}}$ for each algorithm in four experiments, $\pm$ a rough estimate of uncertainty of the value reported (\latin{i.e.} the empirical standard deviation divided by $\sqrt{N}$). In each column, the more accurate algorithms (taking the uncertainty into account) are bolded.
\label{tab.tout}}
\begin{center}
\begin{tabular}{p{0.18\textwidth}@{\hspace{0.025\textwidth}}p{0.16\textwidth}@{\hspace{0.025\textwidth}}p{0.16\textwidth}@{\hspace{0.025\textwidth}}p{0.18\textwidth}@{\hspace{0.025\textwidth}}p{0.205\textwidth}}
\hline\noalign{\smallskip}
Experiment  & S1      & S2        & HSd1            & HSd2 \\
\noalign{\smallskip}
\hline
\noalign{\smallskip}
$\bayes$          & $\sin$  & $\sin$    & HeaviSine       & HeaviSine \\
$\sigma(x)$       & 1       & $x$       & 1               & $x$ \\
$n$ (data) 
				          & 200     & 200       & 2048            & 2048 \\
$\M_n$            & regular & 2 bin sizes & dyadic, regular & dyadic, 2 bin sizes\\
\noalign{\smallskip}
\hline
\noalign{\smallskip}
Mal     & $1.928 \pm 0.04$ & $3.864 \pm 0.02$ & $1.606\pm 0.015$& $1.487\pm 0.011$ \\
Mal+    & $\meil{1.800 \pm 0.03}$ & $4.047 \pm 0.02$ & $1.606\pm 0.015$& $1.487\pm 0.011$ \\
$2-$FCV & $2.078 \pm 0.04$ & $2.542 \pm 0.05$ & $\meil{1.002\pm 0.003}$& $1.184\pm 0.004$ \\
$5-$FCV & $2.137 \pm 0.04$ & $2.582 \pm 0.06$ & $1.014\pm 0.003$& $1.115\pm 0.005$ \\
$10-$FCV& $2.097 \pm 0.05$ & $2.603 \pm 0.06$ & $1.021\pm 0.003$& $1.109\pm 0.004$ \\
$20-$FCV& $2.088 \pm 0.04$ & $2.578 \pm 0.06$ & $1.029\pm 0.004$& $1.105\pm 0.004$ \\
\noalign{\smallskip}
\hline
\noalign{\smallskip}
penEfr     & $2.597\pm 0.07$ & $3.152\pm 0.07$ & $1.067\pm 0.005$& $1.114\pm 0.005$\\
penRad     & $1.973\pm 0.04$ & $2.485\pm 0.06$ & $1.018\pm 0.003$& $\meil{1.102\pm 0.004}$\\
penRHO     & $1.982\pm 0.04$ & $2.502\pm 0.06$ & $1.018\pm 0.003$& $1.103\pm 0.004$\\
penLOO     & $2.080\pm 0.05$ & $2.593\pm 0.06$ & $1.034\pm 0.004$& $1.105\pm 0.004$\\
pen$ 2-$FCV& $2.578\pm 0.06$ & $3.061\pm 0.07$ & $1.038\pm 0.004$& $1.103\pm 0.005$\\
pen$ 5-$FCV& $2.219\pm 0.05$ & $2.750\pm 0.06$ & $1.037\pm 0.004$& $1.104\pm 0.004$\\
pen$10-$FCV& $2.121\pm 0.05$ & $2.653\pm 0.06$ & $1.034\pm 0.004$& $1.104\pm 0.004$\\
pen$20-$FCV& $2.085\pm 0.04$ & $2.639\pm 0.06$ & $1.034\pm 0.004$& $1.105\pm 0.004$\\
\noalign{\smallskip}
\hline
\noalign{\smallskip}
penEfr+     & $2.016\pm 0.05$ & $2.605\pm 0.06$& $1.011\pm 0.003$& $\meil{1.097\pm 0.004}$\\
penRad+     & $\meil{1.799\pm 0.03}$ & $\meil{2.137\pm 0.05}$& $\meil{1.002\pm 0.003}$& $\meil{1.095\pm 0.004}$\\
penRHO+     & $\meil{1.798\pm 0.03}$ & $\meil{2.142\pm 0.05}$& $\meil{1.002\pm 0.003}$& $\meil{1.095\pm 0.004}$\\
penLOO+     & $\meil{1.844\pm 0.03}$ & $\meil{2.215\pm 0.05}$& $\meil{1.004\pm 0.003}$& $\meil{1.096\pm 0.004}$\\
pen$ 2-$FCV+& $2.175\pm 0.05$ & $2.748\pm 0.06$& $1.011\pm 0.003$& $1.106\pm 0.004$\\
pen$ 5-$FCV+& $1.913\pm 0.03$ & $2.378\pm 0.05$& $\meil{1.006\pm 0.003}$& $\meil{1.102\pm 0.004}$\\
pen$10-$FCV+& $1.872\pm 0.03$ & $2.285\pm 0.05$& $\meil{1.005\pm 0.003}$& $\meil{1.098\pm 0.004}$\\
pen$20-$FCV+& $1.898\pm 0.04$ & $2.254\pm 0.05$& $\meil{1.004\pm 0.003}$& $\meil{1.098\pm 0.004}$\\
\hline
\end{tabular}
\end{center}
\end{table}

\subsection*{On the tuning parameters}
The above simulations confirm that the best weights (for accuracy) are Random hold-out ($n/2$) and Rademacher, whereas $V$-fold or leave-one-out weights may be of interest for computational purposes. The second tuning parameter, $\seuilminphl$, may be taken equal to 2 (its ``minimal'' value because terms of the penalty with $n\phl=1$ would be zero) without serious consequences on $C_{\mathrm{or}}$ in practice.

On the contrary, the constant $C \geq C_W$ is quite important, and the best ratio $C/C_W$ strongly depends on $n$, $\sigma$, $\bayes$ and $\M_n$. Moreover, there is no reason for $C_W$(histograms) to be the right non-asymptotic constant in the general algorithm~\ref{def.proc.gal}. Our suggest is to choose $C$ with the so-called ``slope heuristics'', proposed by Birg\'e and Massart \cite{Bir_Mas:2006} for penalties linear in dimension. Their claim is that the optimal penalty is twice the minimal penalty, \latin{i.e.} the one under which the selected model is obviously too large. This leads to estimating the shape of $\penid$ by resampling, and the constant $C$ with the slope heuristics, as follows.
\begin{proc}[Resampling penalization with slope heuristics] \label{def.proc.gal.pente}
\begin{enumerate}
\item Choose a resampling scheme, \latin{i.e.} the law of a weight vector $W$.
\item Compute the following resampling penalty for each $\mM_n$ :
\begin{equation*} 
\pen_0(m) = \Es \croch{  P_n \gamma \paren{ \ERM_m \paren{ \Pnb}} - \Pnb \gamma \paren{\ERM_m \paren{\Pnb}} } \enspace .
\end{equation*}
\item Compute the selected model $\mh(C)$ as a function of $C>0$
\begin{equation*} 
\mh(C) \in \arg\min_{\mM_n} \left\{ P_n \gamma(\ERM_m(P_n)) + C\pen_0(m) \right\} \enspace .
\end{equation*}
\item Choose the minimal $C=\widehat{C}$ such that $D_{\mh(C)}$ is ``reasonably small'', and take $\mh = \mh(2\widehat{C})$.
\end{enumerate}\end{proc}
Step~4 may need to artificially introduce huge models in $\M_n$, all the other ones being considered as ``reasonably small''.
Finally, notice that $C \mapsto \mh(C)$ is piecewise constant with at most $\card(\M_n)$ jumps, so that steps~3--4 have a complexity $\grandO(\card(\M_n))$. As a consequence, the $V$-fold algorithm~\ref{def.proc.gal.pente} is fastly computable.

\section*{Acknowledgements}
I gratefully thank Pascal Massart for many fruitful discussions.
\bibliographystyle{plain}
\bibliography{colt_court}

\end{document}